\documentclass[12pt]{amsart}

\usepackage{latexsym}
\usepackage[all]{xy}
\usepackage{amsmath, amsthm }
\usepackage{amssymb}
\usepackage{amsfonts}
\usepackage{color}
\usepackage{mathtools}
\usepackage[T1]{fontenc}
\usepackage{mathtools}
\usepackage{stmaryrd}
\usepackage{hyperref}
\usepackage{xcolor}
\hypersetup{colorlinks=true,urlcolor=blue,linkcolor=blue,citecolor=blue}

\usepackage{cite}

\newcommand {\Top} {\mathsf{Top}}

\newcommand {\rg} {\mathbb{R}\Gamma}

\newcommand{\ZZ}{\mathbb{Z}}

\newcommand{\Gm}{\mathbb{G}_{m}}
\newcommand{\Ga}{\mathbb{G}_{a}}

\newcommand{\D}{\mathcal{D}}

\newcommand {\Map} {\mathbf{Map}}
\newcommand {\uMap} {\underline{\mathbf{Map}}}

\newcommand {\OO} {\mathcal{O}}

\newcommand {\Spec} {\mathbf{Spec}}

\newcommand  {\dg}     {\mathbf{dg}}

\newcommand  {\St}     {\mathbf{St}}
\newcommand  {\Aff}     {\mathbf{Aff}}

\newcommand  {\ChAff}     {\mathbf{ChAff}}

\newcommand  {\Ab}     {\mathbf{Ab}}
\newcommand  {\AbSch}     {\mathbf{AbSch}}

\newcommand  {\AbChAff}     {\mathbf{AbChAff}}
\newcommand  {\cAbChAff}     {\widehat{\mathbf{AbChAff}}}

\newcommand  {\AbSt}     {\mathbf{AbSt}}

\newcommand  {\dSt}   {\mathbf{dSt}}

\newcommand{\s}{\infty}
\newcommand{\HH}{\mathbb{H}}

\newcommand{\WW}{\mathbb{W}}

\theoremstyle{plain}
\newtheorem{thm}{Theorem}[section]
\newtheorem{df}[thm]{Definition}
\newtheorem{prop}[thm]{Propositon}
\newtheorem{rmk}[thm]{Remark}
\newtheorem{cor}[thm]{Corollary}

\newtheorem{lem}[thm]{Lemma}

\begin{document}

\title{Generalized local jacobians and commutative group stacks}
\author{Bertrand To\"en}

\date{August 2023}

\maketitle

\begin{abstract} 
In \cite[page 109]{MR1942134} Grothendieck sketches the construction of a complex
$J_*(X)$ made of commutative pro-algebraic groups, associated to a smooth variety $X$, and for
which each $J_i(X)$ is a product of local factors called the \emph{local generalized jacobians}. 
The complex $J_*(X)$
has the property to compute Zariski cohomology of $X$ with coefficients in algebraic groups and plays the role of 
a homology of $X$.
The purpose of this note is to recast this construction in the setting of higher algebraic
group stacks for the fppf topology instead. 
For this, we introduce a notion of algebraic homology associated 
to a scheme over a perfect field which is a universal object computing fppf cohomology with coefficients in group schemes. 
We endow this algebraic homology with a filtration by dimension of supports, and prove that,
when $X$ is smooth of dimension $n$, it is $n$-connective with respect to Beilinson's t-structure on the category of 
filtered objects, and that $J_*(X)$ appears as its $n$-th truncation.
In a final
part we partially extend our constructions and results over arbitrary bases.
\end{abstract}

\tableofcontents

\section*{Introduction}

The main purpose of this note is to introduce a notion of \emph{algebraic homology}, 
associated to schemes over a perfect field $k$, which is a universal object 
computing flat cohomology with coefficients in various commutative group schemes. The idea of algebraic homology 
is certainly not new, and was already sketched in a slightly different setting in a letter of Grothendieck to Serre (see letter from 8.9.1960 in \cite[page 109]{MR1942134}), as well as in \cite{grotharchives} (cote 134-2 page 60) 
and in \cite{MR0419452} for instance. However, in the present work we realize algebraic homology 
as a universal (pro)object in the $\s$-category of commutative groups in algebraic higher stacks, whereas in 
the above references it is constructed as a complex of commutative pro-algebraic groups. The main result of this work 
is a comparison statement for smooth schemes, and states that our notion of algebraic
homology recovers Grothendieck's original construction. This comparison is not a formal statement as it relates 
two objects living in two non-equivalent categories, and requires to use an extra piece of data in the form 
of a filtration by dimension of supports. Grothendieck's original construction is shown 
to appear as the $n$-th truncation of algebraic homology, for the Beilinson's t-structure on filtered objects (see for instance \cite{tstructureBeilinson}). 
We also show that both objects can be used in order to compute fppf cohomology with coefficients in 
commutative group schemes (with some affineness and unipotency conditions for higher cohomology groups). This last
statement is true by definition for algebraic homology but requires a proof in the case of Grothendieck's construction 
which is done by gluing together the so-called \emph{local generalized jacobians} appearing 
in our title.

We work over a perfect field $k$, and in the $\s$-topos of 
hypercomplete fppf stacks $\St_k$.
In a first part of this work we introduce the notions of 
\emph{commutative algebraic group stacks over $k$} (see Definition \ref{d1}), which are the commutative group objects
inside the $\s$-category of Artin stacks locally of finite type over $k$ and whose
diagonal is quasi-compact (in a strong sense). These are higher stacks analogues of \emph{locally 
algebraic groups over $k$}, that is commutative group schemes locally of finite type.
The commutative algebraic group stacks can also be characterized as connective complexes of sheaves of abelian 
groups $E$ of the big fppf site, whose cohomology sheaves $H^i(E)$ are all representable by 
locally algebraic groups, and are moreover quasi-projective if $i<0$.
The algebraic homology of a scheme $X$ over $k$ (or any stack over $k$) is then 
defined as follows (see Definition \ref{d4}).

\begin{df}
For a scheme $X$ over $k$, its algebraic homology $\HH^{alg}(X)$ is the pro-object in 
commutative algebraic group stacks which pro-represents the $\s$-functor
$E \mapsto \Map(X,E)$. 
\end{df}

The existence of the pro-object $\HH^{alg}(X)$ in the definition above follows formally from the left exactness
of the $\s$-functor $E \mapsto \Map(X,E)$. However, its existence already carries interesting cohomological information.
For instance, the universal morphism of stacks $X \to \HH^{alg}(X)$  can be seen to 
be a generalization of the Albanese morphism of $X$ which includes higher cohomological 
data in the picture. In a way, $\HH^{alg}(X)$ could also be considered as a derived version of the
classical Albanese variety of $X$ (see end of \S 2).

In order to express $\HH^{alg}(X)$ in more explicit terms we introduce the standard filtration by dimension 
of supports 
$$0 \to F_{dim X}\HH^{alg}(X) \dots \to \dots F_i\HH^{alg}(X) \to F_{i-1}\HH^{alg}(X) \dots \to F_0\HH^{alg}(X)=\HH^{alg}(X)$$ 
where
$F_i\HH^{alg}(X)$ is the part of the algebraic homology which is supported in dimension $i$ in $X$ (see definition \ref{d5}).
The $d$-th graded piece of this filtration is given by the product of all local algebraic homologies
over all the points of dimension $d$ in $X$ (see Proposition \ref{p4})
$$Gr_d\HH^{alg}(X) \simeq \prod_{x\in X^{(d)}}\HH^{alg}_x(X_x).$$
When $X$ is furthermore smooth over $k$, purity implies that 
each graded piece $Gr_d\HH^{alg}(X)$ is in fact $(Dim X-d)$-connective (see Corollary \ref{cp5}). In other words, 
the filtered object $F_{*}\HH^{alg}(X)$ is $0$-connective with respect to Beilinson's t-structure on
filtered objects. Being connective, we get 
a canonical projection of filtered objects 
$$F_{*}\HH^{alg}(X) \longrightarrow H^0_B(F_{*}\HH^{alg}(X))$$
where $H^0_B$ refers to the $0$-th cohomology object for the Beilinson t-structure. The object
$H^0_B(F_{*}\HH^{alg}(X))$ is a genuine complex made of commutative pro-algebraic groups, and as such can be compared to 
Grothendieck's original construction $J_*(X)$. In the complex $J_*(X)$, the degree $d$ part is the 
product of local generalized jacobian $J_x(X)$ over all points $x$ of codimension $d$ in $X$ which are
defined to pro-represent local cohomology in degree $d$
(see Definition \ref{d6}).

Our main comparison result is the following theorem (see Corollary \ref{c2}).

\begin{thm}
If $X$ is smooth over $k$, we have a natural isomorphism of complexes of pro-algebraic groups over $k$
$$H^0_B(F_{*}\HH^{alg}(X)) \simeq J_*(X).$$
\end{thm}

The theorem below, shows that $J_*(X)$ appears as the $0$-th truncation of $F_{*}\HH^{alg}(X)$. However, the natural
projection $F_{*}\HH^{alg}(X) \to J_*(X)$ is not an equivalence in general, and it is thus surprising
that both objects still compute the same cohomologies with unipotent coefficients.

\begin{thm}
Let $i\geq 0$ and $G$ a commutative algebraic group $G$, 
affine if $i>0$ and unipotent if $i>1$. Then, 
the canonical projection $\HH^{alg}(X) \to J_*(X)$ induces isomorphisms
$$Ext^i_{naive}(J_*(X),G) \simeq H^i_{fppf}(X,G).$$
\end{thm}

In the theorem above the $Ext^i_{naive}$ refers to the Yoneda ext-groups computed inside the
derived category of complexes of commutative pro-algebraic groups. These ext-groups differ from the
ext-groups computed as fppf sheaves, and thus are not given by ext-groups computed inside our category 
of commutative group stacks. We note that the theorem above is already stated in
\cite[page 109]{MR1942134} where the fppf topology is replaced with the Zariski topology.

All the results mentioned above are true over a perfect field $k$. We do not know a way
to extend our results and constructions 
over more general base schemes, at least in an interesting manner. The reader will notice that a key 
ingredient for our results is purity for local cohomology of regular schemes with coefficients in 
group schemes. Very general purity theorems have been proven recently in \cite{MR4681144}, but we would rather need 
purity results with coefficients in commutative algebraic group stacks, and we do not see a simple
manner to obtain those from the case of flat group schemes. We hope to be able to come back to this question 
in a future work.

However, when 
restricted to unipotent coefficients, we propose a generalization of algebraic homology 
over any base ring, based on our notion of affine stacks (see \cite{chaff}) and their
abelian counter-parts which we call \emph{affine homology types}. We prove that any 
(small) stack $X$, in particular any scheme, maps to a universal affine homology type called its
affine homology and denoted by $\HH^u(X)$ (see Proposition \ref{p8}
and Definition \ref{d8}). Over a field, $\HH^u(X)$ is the unipotent completion 
of algebraic homology, and its advantage is that it makes sense over any base. Affine homology comes
equipped with a filtration by dimension of supports whose associated graded pieces are given by local 
affine homologies. Moreover, we prove that $\HH^u(X)$ is completely characterized 
by the complex $\mathbb{R}\Gamma(X,\Ga)$ considered as a dg-algebra over $\mathbb{R}End(\Ga)$, which we consider
as a form of Dieudonné theory for affine homology types. \\

To finish this introduction we would like to mention several related works. As we already mentioned
the notion of algebraic homology already appears in \cite[page 109]{MR1942134}) and \cite{MR0419452}.
The results of the present work do not pretend
to be very original and the main point of this note was rather to recast these works in 
the more modern language of higher stacks. The only exception is possibly our last section 
on affine homology, which we think can be of independent interest. 
After posting a first version of this work O. Gwilliam pointed to us the work \cite{MR3697895}
in which a general setting for algebraic homology based on higher commutative group stacks is presented.

In another direction, the results of this work have some limitations. First of all the fact that 
we work over a base field is restrictive. Over more general bases, the $1$-truncation of algebraic
homology essentially appears in the work of S. Brochard in \cite{MR4206615}. Extending the results
of \cite{MR4206615} seems however difficult because of the presence of complicated higher ext groups
between commutative group schemes which prevent to a perfect duality by mapping to $\Gm$.
Also, the approach in \cite{MR4206615} is based on duality, which is an aspect we haven't touched at all 
in this note. Grothendieck already mentions duality in \cite[134-2]{grotharchives}, and 
in degree $0$ this duality of local 
generalized jacobians has been studied in \cite{MR1106897} (it also appears
in \cite{MR3697895}). 

Our comparison requires $X$ to be
smooth, and in his letter Grothendieck already suggests an extension of the definition of $J_*(X)$ to the
singular setting as well, but we could not reconstruct from it anything meaningful at the moment. 

Finally, there are also relations with 
algebraic homotopy types, as been studied for instance in \cite{chaff,kpt,MR3311586,Mondal}. Our notion of affine homology 
is the abelian analogue of these, and the two should be related by means of the
Hurewitz morphism. As a final comment, affine stacks are duals to cosimplicial commutative algebras, 
also called $LSym$-algebras. These are Koszul dual to partition Lie algebras, thanks to \cite{brantner2023deformation}, 
and therefore our notion of affine homology should be closely related to that of
abelian partition Lie algebras. \\

\textbf{Acknowledgments:} In a first version of this work algebraic homology was only studied with some extra 
artificial unipotency conditions. I warmly thank B. Kahn for his comments during a lecture on the content of this work which has
convinced me that these extra conditions were not necessary at all. I am also very grateful to the anonymous referee 
for his very detailed and very useful comments.

\section{Commutative algebraic group stacks}

We work over a perfect base field $k$. 
We denote by $\Aff_k$ the category of affine $k$-schemes, which 
will be endowed with the fppf topology. 

\subsection{Reminders on abelian group schemes}

We denote by $\AbSch_k$, the category of commutative group schemes of finite type over $k$ (or equivalently
quasi-projective over $k$, see for instance \cite[tag 0BF6]{stackproj}). 
The Yoneda embedding provides a fully faithful functor
$$\AbSch_k \hookrightarrow \Ab_k,$$
from $\AbSch_k$ to the category $\Ab_k$ of sheaves of abelian groups on $\Aff_k$ endowed with the $fppf$-topology.
It is well known that this full embedding commutes with the formations of  
kernels and cokernels, and that its essential image is also stable by 
extensions. 

\begin{lem}
The full sub-category of $\Ab_k$ formed by sheaves which are representable
by commutative group schemes of finite type over $k$, is stable by 
taking kernels, cokernels and extensions.
\end{lem}

\textit{Proof.} Using \cite[Exp. $VI_A$ Thm. 5.4.2]{sga3}, and the discussion
after \cite[Exp. $VI_A$ Prop. 5.4.1]{sga3}, we see that the Yoneda embedding is an exact 
functor. Finally, if a sheaf of abelian groups
is an extension of algebraic groups, then it is, by descent, representable by a quasi-separated
algebraic space, and we can conclude by \cite[Lem. 4.2]{MR0260746}. \hfill $\Box$ \\

The Barsotti-Chevalley decomposition theorem (see \cite{Milne-Chevalley}) 
states that any object $G \in \AbSch_k$ fits in an exact sequence
$$0 \to G_0 \to G \to A \to 0,$$ 
where $G_0$ is affine and $A$ is an abelian variety.
Moreover, by \cite[II-5 Cor. 2.3]{MR0302656} $G_0$ itself sits in an exact sequence
$0 \to K \to G_0\to H \to 0,$ 
where $K$ is smooth and affine, and $H$ is 
finite infinitesimal (i.e. $H_{red}=0$).
We also recall from \cite[IV-2 Prop. 2.5]{MR0302656} that $G \in \AbSch_k$ is unipotent if it admits 
a finite filtration by sub-group schemes 
$$\dots\subset G_{i+1} \subset G_i\dots \subset G_{1} \subset G_0=G$$
such that each $G_i/G_{i+1}$ can be realized as a sub-group scheme 
of the additive group $\Ga$. We can moreover arrange things so that 
each $G_i/G_{i+1}$ is a twisted form of $\Ga$, $(\ZZ/p)^r$ for some $r$, or of $\alpha_p$. 
Unipotent group schemes are automatically affine, 
and form a full sub-category of $\AbSch_k$ stable by taking kernels, cokernels and
extensions. 

Finally, an affine commutative group scheme $G$ of finite type over $k$ is always isomorphic to a 
product $G^u\times G^m$ of a unipotent group scheme $G^u$ and a group scheme of multiplicative
type $G^m$. When $k$ is algebraically closed, such a group scheme $G^m$ is always a product 
of multiplicative groups $\Gm$ and roots of unity $\mu_n$. 

\subsection{Reminders on purity for flat cohomology}

We remind here the purity theorem we will use in a crucial manner in order to compare
algebraic homology with Grothendieck's original definition. The most general 
results on purity can be found in \cite{MR4681144}, but for smooth schemes $X$ over a perfect $k$
these can be proven with simple technology, which we recall below. 

Recall that for a point $x$ in a scheme $X$ we define
complex of local cohomology of $X$ at $x$, noted $\mathbb{R}\Gamma_x(X,A)$, with
coefficient in an abelian fppf-sheaf $A \in Ab_k$, as the homotopy fiber of the restriction map
$$\mathbb{R}\Gamma_{fppf}(X_x,A) \to \mathbb{R}\Gamma_{fppf}(X_x-\{x\},A)$$
where $X_x = \Spec\, \OO_{X,x}$ is the local scheme at $x$. It cohomology groups are denoted as usual
by $H^i_x(X,A):=H^i(\mathbb{R}\Gamma_x(X,A))$.

\begin{prop}\label{p0}
Let $X$ be a smooth and connected scheme of dimension $n$ over $k$. For any point $x \in X$ of dimension $d$, 
and any commutative algebraic group scheme $G$ of finite presentation over $k$, we have
$$H^i_x(X,G)=0 \qquad for \, i<n-d.$$
\end{prop}

\textit{Proof.} By Galois descent we can assume that $k$ is algebraically closed.
Let us fix $x \in X$ of dimension $d$ as in the statement. We start by 
the easy observation that the class of groups $G$ satisfying the conclusion of the proposition 
is stable by extensions, and by taking kernels of epimorphic morphisms. We can thus apply 
several standard devissages on $G$ and reduce the statement to the following elementary cases:

\begin{enumerate}
    \item $G$ is unipotent
    
    \item $G$ is an abelian variety
    
    \item $G$ is of multiplicative type.
    
\end{enumerate}

For the case $(1)$ we can use again devissage and thus reduce to the case where $G= \Ga$ is the additive group. 
Similarly, for the case $(3)$ we reduce to the case of the multiplicative group $G=\Gm$. Moreover, 
the case $G=\Ga$ follows from  the usual local coherent duality: the complex $\mathbb{R}\Gamma_x(X,\Ga)$ is cohomologically 
concentrated in the single degree $d$. \\

\textbf{Case $G=\Gm$.} For $d=n$ there is nothing to prove, and when $d=n-1$ the statement simply follows from the fact that
$X_x-\{x\}$ is dense in $X_x$. We thus assume that $d<n-1$.
The case $i=1$ is true by Hartogs, as the restriction map $Pic(X_x) \to Pic(X_x-\{x\})$
is an equivalence of groupoids. For $i=2$ the proposition follows from 
the fact that $H^2_{et}(X_x,\Gm) \to H^2_{et}(X_x-\{x\},\Gm)$ is always injective (see \cite[Cor. 1.10]{MR0244270}).
When $i>2$, the fact that $\Gm$ has torsion cohomology in degree $2$ or higher (see \cite[Prop. 1.4]{MR0244270}) implies that 
$H^i_x(X,\Gm)$ is always torsion. We are thus led to prove the statement for the group $\mu_m$ for some $m$, and by 
devissage we can even assume that $m$ is prime. If $m$ is invertible in $k$, then the purity theorem 
in \'etale cohomology (see for instance \cite[\S VI Thm. 5.1]{MR0559531}) implies the result. When $m=p$ is the characteristic
of the base field, we use the following standard "resolution" of $\mu_p$: on the small fppf site of $X$, there is 
an equivalence
$$\mu_p[1] \simeq (\Omega_{X/k}^{1,cl} \to \Omega_{X/k}^1),$$
where $\Omega_{X/k}^{1,cl} \subset \Omega_{X/k}^{1}$ is the subsheaf of closed $1$-forms and sits in cohomological degrees $0$ (see \cite[\S III Prop. 4.14]{MR0559531}). 
The differential $\Omega_{X/k}^{1,cl} \to \Omega_{X/k}^1$ is here $i-C$, where $i$ is the standard inclusion and $C$ is the 
composition of the projection $\Omega_{X/k}^{1,cl} \to H^1_{DR}(X/k)$ and the Cartier isomorphism
$C : H^1_{DR}(X/k)\simeq \Omega^1_{X/k}$. Both abelian sheaves $\Omega_{X/k}^{1,cl}$ and $\Omega_{X/k}^{1}$
are quasi-coherent (the first one on the Frobenius twists $X^{(1)}$ and the second one on $X$), and thus
satisfy purity. We deduce from this that $H^i_x(X,\mu_p[1]) \simeq H^{i+1}_x(X,\mu_p)\simeq 0$ as soon
as $i < n-d$ and in particular that $H^{i}_x(X,\mu_p) \simeq 0$ when $i<n-d$. \\

\textbf{Case $G$ is an abelian variety.} We first use that all the groups $H^i_{et}(X,G)\simeq H^i_{fppf}(X,G)$ are torsion for $i>0$. Indeed, 
let $i : \eta \to X$ be the generic point of $X$, and $E = \mathbb{R}i_*(i^*(G))$. Then $H^i(X,E)\simeq H^i(\eta,G)$
is Galois cohomology of the field of fractions of $X$ with coefficients in $G$, and thus
is torsion for all $i>0$. Moreover, the Weil extension theorem implies that 
the sheaf $\underline{H}^0(E)$ is canonically isomorphic to $G$ (via the adjunction map $G \to E \simeq \mathbb{R}i_*i^*(G)$). 
We consider the triangle $G[0] \to E \to E_{>0}$, and we notice that all the cohomology sheaves of $E_{>0}$
are torsion sheaves because are also described in terms of Galois cohomology. The fibration
$$\rg(X,G) \to \rg(X,E) \to \rg(X,E_{>0})$$
then implies that all the cohomology groups $H^i(X,G)$ must also be torsion.

When $i=0$, $H^0_x(X,G)=0$ by the density of $X_x - \{x\}$ in $X_x$. Moreover, when $d<n$, 
the restriction $G(X_x) \to G(X_x-\{x\})$  is bijective by Weil's extension theorem, and thus
$H^1_x(X,G)$ identifies with the kernel of $H^1(X_x,G) \to H^1(X_x-\{x\},G)$ and are thus torsion groups.
Finally, by what we have seen $H^i_x(X,G)$ are torsions for $i>1$ as well. Therefore, to prove the 
proposition we can replace $G$ by the finite group scheme $G[n]$ which is the kernel of the multiplication
by $n$ on $G$. This finite group scheme can be split into a unipotent part and a reductive part and thus
have already been treated in the cases before. \hfill $\Box$\\

\subsection{Stacks}\label{sec:stacks}

We denote by $\St_k$ the $\s$-category of (hypercomplete) stacks over the site $(\Aff_k,fppf)$, 
of affine $k$-schemes endowed with the finitely presented flat topology. 
This $\s$-category can be described explicitly as being the $\s$-category of cofibrant and fibrant 
objects in the model category $\Aff_k^{\sim,fppf}$ of presheaves of simplicial sets on $\Aff_k$, endowed with the
local projective model category structure (see for instance \cite[\S 1.1]{chaff} and \cite[\S 2.1]{toenems}).
The objects of $\Aff_k^{\sim,fppf}$ are functors
$$\Aff_k^{op} \to sSet,$$
from the category of affine $k$-schemes to the category of simplicial sets (we neglect size issues here, 
see \cite{chaff} for details concerning universes in this context). For such an object $F$, $X \in \Aff_k$, $i>0$ 
and $s\in F(X)_0$ a $0$-simplex, we have homotopy sheaves $\pi_i(F,s)$, which are sheaves on the site
$\Aff_k/X$ defined as the $fppf$-sheaf associated to the presheaf $(u : Y\to X) \mapsto \pi_i(F(Y),u^*(s))$.
Similarly, we have a sheaf of connected components $\pi_0(F)$, as being associated to the presheaf 
$X \mapsto \pi_0(F(X))$.

The weak equivalences in the model category $\Aff_k^{\sim,fppf}$ are the
morphisms $F \to G$ inducing isomorphism on all possible homotopy sheaves. The $\s$-category $\St_k$ can then be 
realized as the $\s$-category obtained from the genuine category $\Aff_k^{\sim,fppf}$ by formally inverting 
the weak equivalences. A concrete model, up to a natural equivalence of $\s$-categories, consists of the
simplicially enriched category of fibrant and cofibrant objects in $\Aff_k^{\sim,fppf}$. Note moreover, that 
the fibrant objects are precisely the presheaves $F : \Aff_k \to sSet$ satisfying the three conditions below.

\begin{enumerate}
    \item For all $X \in \Aff_k$ the simplicial set $F(X)$ satisfies the Kan condition (i.e. is fibrant as a simplicial 
    set).
    \item For any finite coproduct of affine schemes $X = \coprod_i X_i$, the natural morphism
    $F(X) \to \prod_iF(X_i)$ is an equivalence.
    
    \item For any $fppf$-hypercoverings of affine $k$-schemes $U_* \to X$, the natural morphism
    $$F(X) \to \underset{[n] \in \Delta}{holim}F(U_n)$$
    is a weak equivalence of simplicial sets.
\end{enumerate}

The last two conditions above form the \emph{$fppf$-hyperdescent condition} and is the characteristic property
of \emph{$fppf$-stacks}.
We refer the reader to \cite{chaff,hagI,toenems,htt} for general references on stacks. 

We recall the existence of Postnikov decomposition which will be used several times in the sequel. 
For any stack $F \in\St_k$, and any $n\geq 0$, there exists canonical truncations
$$p_n : F \to \tau_{\leq n}(F)$$
which are characterized by the fact that all homotopy sheaves of $\tau_{\leq n}(F)$
vanish in degree strictly higher than $n$, and moreover the morphism $p_n$ induces 
isomorphisms on all homotopy sheaves in degree less or equal to $n$. For any $n>1$, there exists a cartesian square in 
$\St_k/\tau_{\leq 1}(F)$
$$\xymatrix{
\tau_{\leq n}(F) \ar[r] \ar[d] & \tau_{\leq n-1}(F) \ar[d]^{k_n} \\
\bullet \ar[r] & K(\pi_n(F),n+1)
}$$
where $\pi_i(F)$ is the $n$-th homotopy sheaf of $F$, which is a sheaf of abelian 
groups over the $1$-truncation $\tau_{\leq 1}(F)$. Here $K(\pi_n(F),n+1)$ denotes the Eilenberg-MacLane stack
(relative to $\tau_{\leq 1}(F)$) whose unique non-zero homotopy sheaf is $\pi_n(F)$ in degree $n+1$.
The morphism $k_n$ is called the \emph{$n$-th Postnikov invariant of $F$}.
In our context, the stacks considered will always be abelian group stacks, 
and thus the sheaves $\pi_n(F)$ will be constant over $\tau_{\leq 1}(F)$ and simply be considered as $fppf$-sheaves over
$\Aff_k$.

We denote by $\AbSt_k$ the $\s$-category obtained from $\St_k$ by 
$\ZZ$-linearization. It can be described explicitly as the $\s$-category
of $\s$-functors
$$E : \Aff_k^{op} \to \dg^c_\ZZ,$$
from the category of affine $k$-schemes to the category $\dg^c_\ZZ$ of 
connective complexes of abelian groups, satisfying the 
previously mentioned $fppf$ hyper-descent condition. Any object $E \in \Ab_k$ possesses an underlying 
stack $|E| \in \St_k$ obtained by the Dold-Kan construction, explicitly defined by 
$$|E|(X) := \Map_{\dg^c_\ZZ}(\ZZ,E(X))$$
for any $X \in \Aff_k$, where $\Map_{\dg^c_\ZZ}$ denotes the mapping spaces of the 
$\s$-category $\dg^c_\ZZ$.
This defines a forgetful $\s$-functor $\Ab_k \to \dSt_k$, which 
admits a left adjoint sending $F$ to $\ZZ \otimes F$. For 
$E \in \AbSt_k$ we denote by $\pi_i(E)$ the cohomology sheaves $H^{-i}(E)$, which are canonically identified with 
the homotopy sheaves $\pi_i(|E|)$ of the underlying stack $|E|$ pointed at global point $0$.
Finally, for $E \in \AbSt_k$, the cartesian square coming from the Postnikov decompositions admit 
canonical promotion to cartesian square in $\AbSt_k$
$$\xymatrix{
\tau_{\leq n}(E) \ar[r] \ar[d] & \tau_{\leq n-1}(E) \ar[d]^{k_n} \\
\bullet \ar[r] & \pi_n(E)[n+1]
}$$
(note that the underlying stack of $\pi_n(E)[n+1]$ is $K(\pi_n(E),n+1)$).

\subsection{Commutative algebraic group stacks}\label{sec:commutativegroupstacks}

We remind from \cite[\S 2.2.1]{hagII} the notion of \emph{Artin $n$-stack}. 

\begin{df}\label{d1}
A \emph{commutative algebraic group stack} is an object $E \in \AbSt_k$ such that the
underlying stack $|E| \in \St_k$ is an Artin $n$-stack (for some $n$), locally 
of finite presentation, strongly quasi-separated over $k$ and which can be covered, in the
sense of the fppf topology,
by a countable number of affine $k$-schemes.
They form, by definition, a full sub-$\s$-category 
$$\AbSt^{alg}_k \subset \AbSt_k.$$
\end{df}

We recall that \emph{strongly quasi-separated} means that the diagonal map
$E \to E \times E$ is strongly quasi-compact in the sense of \cite{hagII}. This
means that all the higher diagonal maps
$X \to X^{S^n}$ are quasi-compact for $n\geq 0$. Equivalently, the
inertia stack $IE \to E$ is required to be strongly quasi-compact over $E$.
The countability assumption is here to ensure that commutative algebraic group
stacks form an essentially small $\s$-category.

As explained in \cite[Cor. 2.9]{Gerbe} an Artin stack $E$ locally of finite presentation 
which is strongly quasi-separated is such that 
for all field valued point $s \in E(K)$, all its homotopy sheaves $\pi_i(E,s)$ 
are quasi-compact and quasi-separated algebraic spaces of finite ype over $K$, and thus are 
representable by group schemes thanks to \cite[Lem. 4.2]{MR0260746}. 

As we are working over a field, the converse also holds. \\

\begin{prop}\label{p1}
An object $E \in \AbSt_k$ is a commutative algebraic group stack if and only if
each homotopy sheaf $\pi_i(E) \in Ab_k$ is representable by a group scheme
locally of finite type over $k$, which is moreover
quasi-projective for $i>0$ and vanishes for all but a finite number of indices $i$.
\end{prop}

\textit{Proof.} 
Assume first that $E$ is a commutative algebraic group stack. 
In the definition \ref{d1}. We
have to prove that all the homotopy sheaves $\pi_i(E)$
are representable. We consider the inertia stack $IE$ of $E$. We already saw that 
$\pi_i(E)$ are group schemes for $i>0$ by means of \cite[Cor. 2.9]{Gerbe}\footnote{the result of \cite[Cor. 2.9]{Gerbe}
assumes the quasi-compactness of $E$, but the result remains valid 
in our situation as any point of $E$ with value in an affine scheme will factor through a quasi-compact
open substack $E' \subset E$.}, 
it remains to treat the case $i=0$.
As $E$ 
is an abelian group object we have splitting $IE\simeq E \times \Omega_0E$, where
$\Omega_0E$ is the pointed loop stack at $0$. As $k$ is a field, we thus have that 
$IE \to E$ is a fppf cover which is strongly quasi-compact. We can thus 
use \cite[Prop. 2.4]{Gerbe}, which remains valid only assuming the diagonal to 
be strongly quasi-compact, and deduce that $E$ is a gerbe. In particular, 
$\pi_0(E)$ is representable by an algebraic space locally of finite type over $k$.
This is a group object in algebraic space, and is moreover quasi-separated 
by assumptions on $E$, and by \cite[Lem. 4.2]{MR0260746} is representable by 
a group scheme locally of finite type. 

Conversely, assume that $E$ satisfies the conditions of the proposition. 
We will prove by induction on $i$ that the $i$-th homotopical truncation $\tau_{\leq i}E$
is an Artin $i$-stack locally of finite presentation and strongly quasi-separated over $k$.
For $i=0$ this is obvious, as $\tau_{\leq 0}E \simeq \pi_0(E)$ is a scheme locally of finite type by assumption.
Consider $\tau_{\leq i}E \to \tau_{\leq i-1}E$, which fits into a cartesian square
$$\xymatrix{
\tau_{\leq i}E \ar[r] \ar[d] & \tau_{\leq i-1}E \ar[d] \\
 0 \ar[r] & \pi_i(E)[i+1].}$$
By assumption, $\pi_i(E)$ is a group scheme of finite type over $k$ for $i>0$, and thus
$K(\pi_i(E),i+1)$ is an Artin $(i+1)$-stack locally of finite presentation and strongly quasi-separated over $k$.
By induction, this implies that $\tau_{\leq i}E$ is an Artin $i$-stack locally of finite presentation
and strongly quasi-separated.
\hfill $\Box$ \\

By using the long exact sequence in homotopy, we see that the proposition \ref{p1}
easily implies the following corollary.

\begin{cor}\label{cp1}
The full sub-$\s$-category $\AbSt^{alg}_k \subset \AbSt_k$ is  
is stable by finite limits, retracts and the homotopical truncation 
$\s$-functors $\tau_{\leq i}$.
\end{cor}

\begin{rmk}
\emph{
We caution the reader that if $E$ is a commutative algebraic group stack
then $E[1]=BE$ is generally not a commutative algebraic group stack anymore. 
It is so only if $\pi_0(E)$ is moreover assumed to be 
of finite type.}
\end{rmk}

We will also use occasionally the following bigger category.

\begin{df}\label{d1'}
A \emph{commutative locally algebraic group stack} is an object $E \in \AbSt_k$ 
such that for all $i$ the sheaf $\pi_i(E)$ is representable by an algebraic space locally
of finite type. The $\s$-category of commutative locally algebraic group stacks
is denoted by $\AbSt^{lalg}_k \subset \AbSt_k$.
\end{df}

Unlike algebraic group stacks, the locally algebraic stacks are stable
by suspension: if $E \in \AbSt_k^{lalg}$, then $BE$, computed in $\AbSt_k$, remains in 
$\AbSt_k^{lalg}$. More generally, $\AbSt_k^{lalg}$ is stable by finite limits and colimits,
as this can be easily seen using the long exact sequence in homotopy.

\section{Algebraic homology}

In this section we define algebraic homology and study its basic properties. \\

We let $X$ be any stack over $k$. We consider the $\s$-functor
$$\AbSt_k^{alg} \longrightarrow \Top$$
defined by $E \mapsto \Map(X,E)$, and where $\Top$ denotes the $\s$-category of Kan simplicial sets. 
This $\s$-functor will be denoted by 
$E\mapsto \HH(X,E)$. By definition it obviously commutes with finite limits, and thus
defines a pro-object in $\AbSt_k^{alg}$ (see \cite[Def. 7.1.6.1]{htt} and \cite[Rem. 7.1.6.2]{htt})
$$\HH^{alg}(X) \in Pro(\AbSt_k^{alg}).$$

\begin{df}\label{d4}
With the notation above, the object $\HH^{alg}(X)$ is called the
\emph{algebraic homology of $X$ (relative to $k$)}. 
\end{df}

Note that the 
definition \ref{d4} makes sense for any stack $X$ over $k$, possibly not even representable
by an Artin stack. However, we will mainly be interested in the case where
$X$ is a scheme of finite type over $k$. Note also that 
algebraic homology is a relative notion and depends on the base field $k$. It should therefore
be noted by $\HH^{alg}(X/k)$ when mentioning the base field is important.

By definition we have in particular that for any 
commutative group scheme $G$ of finite type over $k$
$$H^i_{fppf}(X,G) \simeq \pi_0(\Map_{\St_k}(X,K(G,i)) \simeq
\pi_0(\Map_{Pro(\AbSt_k)}(\HH^{alg}(X),K(G,i))).$$

\begin{rmk}
\emph{The notion above of algebraic homology is very close to that of flat homology of 
\cite{MR0419452}, but is different. The major difference is that 
these two objects live in two different, non-equivalent, $\s$-categories. Flat homology
lives in the $\s$-category of complexes of commutative pro-algebraic groups. This has 
a canonical $\s$-functor to $\AbSt_k$, by considering the stack  pro-represented by 
a complex of pro-algebraic groups, but this $\s$-functor is known not to be 
fully faithful. This is due to the existence of two different, non-isomorphic, 
definitions of ext-groups of commutative algebraic groups: Yoneda ext-groups inside algebraic groups and
ext-groups of the corresponding fppf-sheaves. It is known that these two notions are different (see for 
instance \cite{MR0255550}).}
\end{rmk}

By definition the algebraic homology possesses the following basic properties. \\

\textbf{Functoriality.} The construction $X \mapsto \HH^{alg}(X)$ is functorial in $X$, 
and defines an $\s$-functor from Artin $k$-stacks to pro-objects in derived commutative algebraic
groups. Moreover, the identity of $\HH^{alg}(X)$ defines a canonical 
class $\alpha_X \in \HH(X,\HH^{alg}(X))$, and thus a canonical morphism of pro-stacks
$$\alpha_X : X  \to \HH^{alg}(X),$$
which is functorial in $X$ as well. \\

\textbf{Homology of the point.} We write $\HH^{alg}(k)$ for $\HH^{alg}(\Spec\, k)$. It is easy to 
see that $\HH^{alg}(k)\simeq \ZZ$ is the constant fppf sheaf associated to $\ZZ$ and pro-constant as a pro-object. \\

\textbf{Homology groups.} The pro-object $\HH^{alg}(X)$ possesses pro-homology
groups $\HH_i^{alg}(X):=\pi_i(\HH^{alg}(X))$, 
which are commutative pro-algebraic groups over $k$. By definition 
of our notion of commutative algebraic group stacks 
$\HH_0^{alg}(X)$ is a pro-group scheme locally of finite presentation, and 
all the $\HH_i^{alg}(X)$ are pro-group schemes of finite type for $i>0$.

In general, the projection $X \to \Spec\, k$ provides a canonical morphism
$\HH^{alg}(X) \to \HH^{alg}(\Spec\, k)\simeq \ZZ$. 
We thus have a canonical fibration sequence of commutative algebraic group stacks
$$\HH^{alg}(X)^{(0)} \to \HH^{alg}(X) \to \HH^{alg}(k)\simeq \ZZ,$$
where $\HH^{alg}(X)^{(0)}$ is by definition the \emph{reduced algebraic homology of $X$}.
In general the fibration sequence above does not split (see below), and a choice of a 
rational point $x \in X(k)$ provides a canonical splitting $\HH^{alg}(X) \simeq \ZZ \oplus 
\HH^{alg}(X)^{(0)}$ and corresponding splitting on homology groups. 

Note also that by the universal property, if $X$ is geometrically connected then the pro-\'etale group of connected components of $\HH_0^{alg}(X)$
must be isomorphic to the constant sheaf $\ZZ$ by the above morphism $\HH_0^{alg}(X) \to \ZZ$. \\

\textbf{Relation to the Albanese variety.}
The pro-algebraic group $\HH_0^{alg}(X)$ is directly related to the Albanese variety of $X$. Indeed, 
for any commutative algebraic group $G$, we have a canonical bijection
$$H^0(X,G)=Hom(X,G) \simeq Hom(\HH_0^{alg}(X),G).$$
Fixing a point $x \in X(k)$, we get a decomposition 
$\HH_0^{alg}(X) \simeq \HH_0^{alg}(X)^{(0)} \oplus \ZZ$. 
When 
$X$ has no $k$-point, the projection $X \to \Spec\, k$ provides a (a priory non-split) short exact 
sequence
$$0 \to \HH_0^{alg}(X)^{(0)} \to \HH_0^{alg}(X) \to \ZZ \to 0.$$
We can consider the universal morphism $X \to \HH^{alg}(X)$, and consider the 
composite morphism $X \to \HH_0^{alg}(X) \to \ZZ$. As $\HH^{alg}(-)$ is functorial, this morphism
is the composition $X \to \Spec\, k \to \HH^{alg}(\Spec\, k)=\ZZ$. By definition the second morphism
$\Spec\, k \to \ZZ$ exhibits $\ZZ$ as the free sheaf of abelian groups generated by 
$\Spec\, k$, and thus corresponds to the global section $1 \in \ZZ(k)$. As a consequence,
the composition $X \to \ZZ$ is the global function $1$ on $X$.

As a consequence of the discussion above, the morphism $X \to \HH^{alg}(X)$ factors through a canonical morphism
$$alb_X : X \to \HH_0^{alg}(X)^{(1)}$$
where $\HH_0^{alg}(X)^{(1)}$ is the fiber at $1$ of the short exact 
sequence above, considered naturally as a torsor over the group $\HH_0^{alg}(X)^{(0)}$.
This morphism is a generalization of the usual Albanese morphism to the Albanese variety.
When $X$ is smooth and proper, $\HH_0^{alg}(X)^{(1)}$ can be identified with 
the usual Albanese variety, and $alb_X$ with the usual the Albanese map, as these two objects
satisfy the same universal property.
In general, the universal morphism $X \to \HH^{alg}(X)$ can be understood as
a higher cohomological generalization of the Albanese morphism as being 
the universal map towards a torsor over an object in $Pro(\AbSt_k^{alg})$. \\

\textbf{Relation to the Albanese stack.}
When $X$ is an Artin stack of finite type over $X$, the first Postnikov truncation
$\HH_{\leq 1}^{alg}(X)$ is an abelian group object in algebraic pro-Artin $1$-stacks. 
This object is closely related to $D(Pic(X))$, the dual of the Picard stack of $X$ (relative to $k$)
studied in \cite{MR4206615}. Indeed, $Pic(X):=\uMap_{St_k}(X,B\Gm)$ is
the internal $Hom$ stack from $X$ to $B\Gm$, and we thus have a canonical evaluation morphism
$ev : X \times Pic(X) \to B\Gm$
which defines a canonical morphism
$$X \to \uMap_{Ab_k}(Pic(X),B\Gm)=D(Pic(X)),$$
where $\uMap_{Ab_k}$ denotes the stack of abelian group stacks morphisms.
When $X$ is proper and the dual $D(Pic(X))$ is representable by an algebraic group stack, there is 
a fibration sequence of commutative group stacks 
$$0 \to D(\pi_0 Pic(X)) \to D(Pic(X)) \to \ZZ \to 0.$$
This fibration sequence is obtained by mapping towards $B\Gm$
the fibration sequence
$$0 \to B\Gm=Pic(\Spec\, k) \to Pic(X) \to \pi_0Pic(X) \to 0.$$
 
The choice of a global point $x \in X(k)$ defines a splitting of these two sequences, and thus
the evaluation map $X \to D(Pic(X))$ can be composed to a morphism $X \to D(\pi_0Pic(X)$. 
By the universal property of algebraic homology we get a canonical morphism of
pro-objects in commutative algebraic $1$-stacks
$$\HH_{\leq 1}^{alg}(X)^{(0)} \to D(\pi_0Pic(X)).$$
In absence of global points we get a morphism $\HH_{\leq 1}^{alg}(X)^{(1)} \to D(\pi_0Pic(X))^{(1)}$
where $\HH_{\leq 1}^{alg}(X)^{(1)}$ (resp. $D(\pi_0Pic(X))^{(1)}$) is the fiber at $1$ of the canonical fibration sequence
$$\HH_{\leq 1}^{alg}(X)^{(0)} \to \HH_{\leq 1}^{alg}(X) \to \HH_{\leq 1}^{alg}(\Spec\, k)$$
$$\left( \textrm{resp.} \; D(\pi_0Pic(X)) \to D(Pic(X)) \to \ZZ \right).$$
This map is an isomorphism when $X$ satisfies some extra conditions
expressed in \cite{MR4206615} (e.g. when $X$ is smooth and proper over $k$), 
as both objects $\HH_{\leq 1}^{alg}(X)^{(1)}$ and $D(\pi_0Pic(X))^{(1)}$ then satisfy the same
universal property. This provides a direct relation between algebraic homology 
and the Albanese morphism constructed in \cite{MR4206615}.\\

\textbf{Change of base fields.} The construction $X \mapsto \HH^{alg}(X/k)$ is compatible with 
change of base field in the following sense. Let $k'/k$ be a finite field extension, and $X':=X \times_k k'$ the base
change of $X$ over $\Spec\, k'$. By universal property we have a canonical morphism of pro-objects
$$\HH^{alg}(X'/k') \to \HH^{alg}(X/k)\times_k k'.$$
We claim that this morphism is in fact an equivalence of pro-objects. Indeed, for any $G \in AbSt_{k'}^{alg}$, 
we have
$$Map_{\St_{k'}}(\HH^{alg}(X/k)\times_k k',G) \simeq Map_{\St_k}(\HH^{alg}(X/k),p_*(G)),$$
where $p_*(G)$ is the direct image (i.e. \emph{restriction à la Weil}) of $G$ along the morphism
$p : \Spec\, k' \to \Spec\, k$. Because $p$ is finite and etale it is easy to see that 
the stack $p_*(G)$ is again a commutative algebraic group stack. Therefore, by the universal property of algebraic homology we get
$$Map_{\AbSt_{k'}}(\HH^{alg}(X/k)\times_k k',G) \simeq Map_{\AbSt_k}(\HH^{alg}(X/k),p_*(G))$$
$$\simeq Map_{\St_k}(X,p_*(G))
\simeq Map_{\St_{k'}}(X'\times_k k',G) \simeq Map_{\AbSt_{k'}}(\HH^{alg}(X'/k'),G),$$
as required.

\section{Grothendieck's generalized local jacobians}

In this section we provide $\HH^{alg}(X)$ with 
a natural filtration, induced by the standard filtration by codimension
on cohomology. The associated graded to this filtration possesses a description in local 
terms involving local cohomology. When  working with smooth schemes over $k$, the local factors will
be related to Grothendieck's original construction of local generalized jacobians
from his letter to Serre (letter from 9.8.1960, \cite{MR1942134} page 109). \\

In this section $X$ is a scheme of finite type over $k$. We denote by $X^{(d)}$ the subset
of $X$ consisting of all points of dimension $d$ in $X$. Similarly, $X^{(\geq d)}$ will be
the subset of points of dimension at least equal to $d$.

\subsection{The filtration by (co)dimension}

We let $d\geq 0$ be an integer, and consider $Op_d(X)$ the poset of
opens $U \subset X$ containing $X^{(\geq d)}$. 
An open $U$ belongs to $Op_{d}(X)$ if and only if the closed complement $X-U$ has dimension 
at most $d-1$. The only element in $Op_0(X)$ is $X$ itself, whereas
$\emptyset \in Op_{dim X+1}(X)$, and obviously $Op_{d}(X) \subset Op_{d+1}(X)$.
We define 
$F_d\HH^{alg}(X)$ to be 
$$F_d\HH^{alg}(X):=\lim_{U \in Op_{d}(X)} \HH^{alg}(U)$$
where the limit is taken over the poset $Op_{d}(X)$ and 
inside the $\s$-category $Pro(\AbSt_k^{alg})$ of pro-objects in 
commutative algebraic group stacks. Note that $F_{0}\HH^{alg}(X)\simeq \HH^{alg}(X)$, 
and as $\emptyset$ is initial in $Op_{dim X+1}(X)$ we have $F_{dim X+1}\HH^{alg}(X)\simeq 0$.

The object $F_d\HH^{alg}(X)$
comes equipped with a canonical morphism $F_{d+1}\HH^{alg}(X) \to F_{d}\HH^{alg}(X)$, 
induced from the inclusion of posets $Op_{d}(X) \to Op_{d+1}(X)$, and we thus 
obtain this way a filtration 
$$0=F_{dim X+1}\HH^{alg}(X)=\HH^{alg}(X) \to \dots \to F_{d}\HH^{alg}(X) \to F_{d-1}\HH^{alg}(X) 
\to \dots \to F_{0}\HH^{alg}(X)=\HH^{alg}(X).$$

\begin{df}\label{d5}
The above filtration is called the  \emph{filtration
by dimension} on $\HH^{alg}(X)$.
\end{df}

The $d$-th graded piece $Gr_d\HH^{alg}(X)$ are defined as the cofiber
of $F_{d+1}\HH^{alg}(X) \to F_{d}\HH^{alg}(X)$, computed in
the stable $\s$-category $Pro(\AbSt_k)$ (see comment \ref{r5}) which is dual to local 
cohomology in dimension $d$ on $X$. More precisely, we have the following
universal characterization of $Gr_d\HH^{alg}(X)$.

\begin{prop}\label{p4}
With the notation above, 
and for any commutative algebraic group
stack $G$, there exists a canonical equivalence, functorial in $G$
$$\Map(Gr_d\HH^{alg}(X),G) \simeq \bigoplus_{x\in X^{(d)}}\HH_x(X,G),$$
where the sum runs over the set $X^{(d)}$ consisting of 
all points of $X$ of dimension $d$, and where
$\HH_x(X,G)$ is the space of local cohomology defined as the fiber of $\Map(X_x,G) \to \Map(X_x-\{x\},G)$
with $X_x=\Spec\, \OO_{X,x}$.
Equivalently, we have a canonical identification
$$Gr_d\HH^{alg}(X) \simeq \prod_{x\in X^{(d)}}\HH_{x}^{alg}(X_x)$$
where $\HH_{x}^{alg}(X_x)$ is local algebraic 
homology at $x$, defined as the cofiber of the natural morphism
$\HH^{alg}(X_x-\{x\}) \to \HH^{alg}(X_x)$.
\end{prop}

\textit{Proof.} For any inclusion of open $V\subset U$, with $U\in Op_d(X)$ 
and $V\in Op_{d+1}(X)$, the closed complement $U-V$ only has a finite number of 
points of dimension $d$. Indeed, if $(U-V)^{(d)}$ is infinite, then $U-V$ must be
of dimension at least $d+1$ and thus must contain at least one point of dimension $d+1$, which 
contradicts $V \in Op_{d+1}(X)$. Note also that for $x\in (U-V)^{(d)}$, 
we have $X_x\cap V = X_x-\{x\}$, 
as $x$ is the only point of $X_x$ of dimension $d$ in $X$.
We thus get this way a commutative square of schemes
$$\xymatrix{
\coprod_{x\in (U-V)^{(d)}} X_x-\{x\} \ar[d] \ar[r] &  V \ar[d] \\
\coprod_{x\in (U-V)^{(d)}} X_x \ar[r] & U,}$$
and thus an induced commutative square on algebraic homologies, which induces a map on the corresponding vertical cofibers
$$\prod_{x\in (U-V)^{(d)}}\HH_{x}^{alg}(X_x) \longrightarrow \HH^{alg}(U)/\HH^{alg}(V).$$
By passing to the limit over the pairs $V\subset U$ with $U \in Op_d(X)$ and $V\in Op_{d+1}(X)$
we find the required morphism
$$\phi : \prod_{x \in X^{(d)}} \HH_{x}^{alg}(X_x) \longrightarrow 
\lim_{V\subset U}\HH^{alg}(U)/\HH^{alg}(V) = Gr_d\HH^{alg}(X).$$

It remains to prove that this morphism $\phi$ is an equivalence. For this, by functoriality in $X$, 
we localize the construction on the small Zariski site of $X$. Both 
constructions
$$U \mapsto \prod_{x\in U^{(d)}}\HH_{x}^{alg}(U_x) \qquad
U \mapsto Gr_d\HH^{alg}(U)$$
are costacks with values in $Pro(\AbSt_k)$ on the small Zariski site $X_{Zar}$. It is therefore enough
to check that the above morphism induces an equivalence on the fibers
at a given point $x \in X$. This fiber is clearly $0$ if $x$ is not of dimension $d$. If $x$
is of dimension $d$ the fiber of the right hand side is equivalent to 
$$\lim_{x\in U \in Op_d(X)}\HH^{alg}(U)/\HH^{alg}(U-\overline{\{x\}}).$$
The fiber at $x$ of the left hand side is clearly $\HH_{x}^{alg}(X_x)$.
We are thus reduced to show that for any commutative algebraic group
stack $G\in \AbSt_k^{alg}$, the canonical morphisms
$$\underset{x\in U}{colim}\Map(U,G) \to \Map(X_x,G) \qquad
\underset{x\in U}{colim}\Map(U-Z_x,G) \to \Map(X_x-\{x\},G)$$
are equivalences. But this is true because, by definition of commutative 
group stacks, $G$ is an Artin stack of finite presentation, and thus
its functor of points sends filtered colimits of commutative rings to 
filtered colimits.
\hfill $\Box$ \\

\begin{rmk}\label{r5}
\emph{As we already saw, the commutative algebraic group stacks are not stable
by taking cofibers. So, a priory $Gr_d\HH^{alg}(X)$  
merely belongs to $Pro(\AbSt_k^{lalg})$, the
$\s$-category of pro-objects in commutative locally algebraic groups stacks. However,
it is easy to see by universal property that for any connected smooth variety $U$ over $k$, 
the group of connected components of $\HH_0^{alg}(U)$ is canonically isomorphic to $\ZZ$, as
any morphism from $U$ to a discrete group $G$ must be constant. As a result, 
for any morphism $U \to V$ of such varieties, the induced morphism $\HH_0^{alg}(U) \to \HH_0^{alg}(V)$
is an isomorphism on the group of connected components. This
easily implies that  
$Gr_d\HH^{alg}(X)$ is a pro-object in $\AbSt^{alg}_k$.}
\end{rmk}

As the objects $Gr_d\HH^{alg}(X)$ are pro-objects in $\AbSt_k^{alg}$, their homotopy sheaves are
pro-sheaves of abelian groups. In the results below we will 
simply say that $Gr_d(\HH^{alg}(X))$ is \emph{$m$-connective} for some $m$, to 
signify that the pro-sheaves $\pi_i(Gr_d\HH^{alg}(X))$ vanish for $i< m$.

\begin{prop}\label{cp5}
Let $X$ be a connected and smooth $k$-scheme of dimension $n$. Then, for all integer $d$
the object $Gr_d\HH^{alg}(X)$ is $(n-d)$-connective: its 
pro-homotopy groups $\pi_i(Gr_d\HH^{alg}(X))$ vanish for $i<n-d$.
\end{prop}

\textit{Proof.} When $d=n$ there is nothing to prove as by construction all $Gr_d\HH^{alg}(X)$ are $0$-connective. For 
$d<n$ we have to show that $\pi_i(Gr_{d}\HH^{alg}(X))=0$ for $i<n-d$. 
Equivalently we
have to prove that any 
$G\in \AbSt_k^{alg}$ which is $(n-d-1)$-truncated, then 
$\Map_{\AbSt_k}(Gr_d\HH^{alg}(X),G) \simeq *$. By the proposition 
\ref{p4} this is equivalent to say that 
$\HH_x(X,G)\simeq 0$ for any point $x\in X$ of dimension $d$. 
We then proceed by 
Postnikov decomposition on $G$, and we thus reduce the corollary 
to the following two statements.
\begin{enumerate}
    \item For all $x\in X^{(d)}$, 
    and all abelian group scheme $G$ locally of finite type, we have $H_x^0(X,G)=0$.
    
    \item For all $x \in X^{(d)}$ and all abelian group scheme $G$ of finite type,
    we have $H^i_x(X,G)=0$ for all $i<n-d$
\end{enumerate}

The first of this statement is obviously true because $d<n$ and thus $X_x-\{x\}$ is dense in $X_x$.
The second statement has already been reminded in the proposition \ref{p0}.
\hfill $\Box$ \\

\subsection{Grothendieck's construction}

We will now compare the general notion of algebraic homology with 
the construction sketched by Grothendieck in \cite[page 109]{MR1942134} 
(which also appears in a letter to L. Breen, see \cite[cote 134-2]{grotharchives}  page 60). As we are working with
the fppf topology, we will have to modify slightly 
the presentation that was originally done by Grothendieck in the context of the Zariski topology.
For this, and all along this section, we will assume that $X$ is
a smooth scheme, connected and of dimension $n$ over $k$. 
Let $x \in X^{(d)}$ be a point of dimension $d$ on $X$.

We consider the functor $G \mapsto H_x^{n-d}(X_x,G)$, (fppf) local cohomology at $x$,
defined on the category of abelian group schemes (separated and of finite type) over $k$. By purity
(Proposition \ref{p0})
we have $H_x^{n-d-1}(X_x,G)\simeq 0$ and thus the above functor is left exact in $G$.
It is therefore pro-representable by an abelian pro-group scheme $J_x(X)$.

For an integer $m$ we define
$$J_m(X):=\prod_{x\in X^{(n-m)}}J_x(X)$$
which is a pro-object in the category of commutative group schemes locally of finite type over $k$. 
We denote by $J_*(X)$ the graded
object $\oplus_m J_m(X)$. It comes equipped with a homological differential
$J_m(X) \to J_{m-1}(X)$ defined as follows. For a fixed $z \in X^{(n-m+1)}$, 
we define a morphism 
$$d_z : \prod_{x\in X^{(n-m)}}J_x(X) \to J_z(X)$$
by defining a natural transformation on the functor they pro-represent. 
The right hand side pro-represents the functor $G \mapsto H^{m-1}_z(X_z,G)$, 
whereas the left hand side pro-represents the functor
$G \mapsto \oplus_{x\in X^{(n-m)}}H^{m}_x(X_x,G)$ (note that 
because $G$ is a pro-constant object we have 
$Hom(\prod_{x\in X^{(n-m)}}J_x(X),G) \simeq \oplus_{x\in X^{(n-m)}}Hom(J_x(X),G)$).

In order to define
$d_z$ we thus have to construct a morphism on local cohomologies
$$H^{m-1}_z(X_z,G) \longrightarrow \oplus_{x\in X^{(n-m)}}H^{m}_x(X_x,G),$$
functorial in $G$. This morphism is the usual differential in the
Cousin complex (see \cite[\S 1]{MR1466971}). Any element $\alpha\in  H^{m-1}_z(X_z,G)$
comes from an element $\alpha_U \in H^{m-1}_Z(U,G)$
for an open $U \subset X$ containing $z$, and where $Z$ is the closure of $z$ in $U$.
We denote by $x_1,\dots,x_q$ the points of $X$ of dimension $n-m$ which are
specializations of $z$ but not belonging to $U$, and
we let $T=X-U$. We consider the boundary map
$H^{m-1}(U,G) \to H^{m-1}_T(X,G)$, together with the restriction morphism 
$H^{m}_T(X,G) \to H^{m-1}_x(X_x,G)$
along
the inclusion of pairs $(X_x,\{x\}) \subset (X,T)$. We get this way a composed morphism
$$H^{m-1}_Z(U,G) \to H^{m-1}(U,G) \to H^{m}_T(X,G) \to \bigoplus_{1\leq i\leq q} 
H^{m}_{x_i}(X_{x_i},G).$$
The image of $\alpha_U$ by the morphism is by definition $d_z(\alpha)$. It is well
known that this morphism is independant of the choice of $\alpha_U$ and 
squares to zero.

\begin{df}[Grothendieck]\label{d6}
Let $X$ be a smooth and connected $k$-scheme of dimension $n$ as above.
\begin{enumerate}
    \item For any $x \in X$ of dimension $n-m$, the pro-algebraic group 
    $J_x(X)$ is called \emph{Grothendieck's local generalized jacobian of $X$ at $x$}.
    \item The pro-algebraic group $J_m(X):=\prod_{x\in X^{(n-m)}}J_x(X)$
    called \emph{Grothendieck's local generalized jacobian of $X$}.
    \item The complex of local generalized jacobians, with the differentials
    defined above, will be denoted by 
    $$J_*(X):=\xymatrix{J_{n}(X) \to J_{n-1}(X) \to \dots \to J_1(X) \to J_0(X)}.$$
\end{enumerate}
\end{df}

As a first relation between the complex $J_*(X)$ and algebraic homology, we have the
following direct observation, which is a direct consequence of the proposition \ref{cp5}
and and the definitions of $J_*(X)$. 

\begin{cor}\label{c2}
Let $X$ be a smooth $k$-scheme of dimension $n$. Then, the first non-zero homotopy groups of 
the graded pieces of $\HH^{alg}(X)$ are the Grothendieck's local generalized jacobians
$$\pi_{d}(Gr_{d}\HH^{alg}(X)) \simeq J_{n-d}(X).$$ 
\end{cor}

By the purity result \ref{cp5}, Grothendieck's local generalized jacobians $J_d(X)$
therefore appear as the lowest non-trivial homotopy group of $Gr_{(n-d)}\HH^{alg}(X)$. 
However, the graded pieces $Gr_{(n-d)}\HH^{alg}(X)$ are in general not concentrated in a single degree, 
and thus $Gr_{(n-d)}\HH^{alg}(X)$ and $J_d(X)$ differ.
We will now see however that the complex $J_*(X)$ is essentially enough to recover the cohomology of $X$
with coefficients in certain commutative algebraic groups
as stated in \cite[page 109]{MR1942134}.

The corollary \ref{c2} can also be phrased by stating that the complex of
pro-algebraic groups $J_*(X)$ is the $n$-th truncation of the $n$-connective filtered object
$\HH^{alg}(X)$ with respect to the Beilinson's t-structure on filtered objects. We
remind that in our context of decreasing filtrations, this t-structure
is characterized by the fact that the connective objects are the filtered
objects $F_*E$ such that $Gr_p(E)$ is $-p$-connective for all $p$: $\pi_i(Gr_p(E))=0$ for all $i<-p$ (see for instance \cite{tstructureBeilinson}).
As a result, we do have a canonical projection of filtered objects
$$q : \HH^{alg}(X) \to J_*(X)$$
where the filtration on the right hand side is the naive filtration 
$F_dJ_*(X):=J_{\geq n-d}(X) \subset J_*(X)$.
This projection is not an equivalence, as the graded pieces $Gr_d\HH^{alg}(X)$ are not
purely concentrated in a single degree in general. However, the projection 
$q$ induces an isomorphism on certain ext groups showing that the complex
$J_*(X)$ can be used in order to compute the individual cohomology groups
$H^i_{fppf}(X,G)$ with coefficients in a commutative algebraic group $G$ (affine if $i>0$, and
unipotent if $i>1$). In order to explain this, we need to first introduce
a comparison functor, between the naive derived category of commutative group schemes
and our category of commutative derived group stacks. 

We know that the category $\AbSch_k$, of commutative group schemes 
of finite type over $k$ is an abelian category. Its category of pro-objects
$Pro(\AbSch_k)$ is known to be abelian and to possess enough projective objects (see \cite{MR0266929}).
As such it possesses a bounded derived category
$\D^b(Pro(\AbSch_k))$, obtained by localizing the category of bounded complexes in $Pro(\AbSch_k)$
with respect to the quasi-isomorphisms. Because the functor sending 
$G \in Pro(\AbSch_k)$ to the 
fppf-sheaf it represents is exact, there is a canonical $\s$-functor
$$h : \D^b(Pro(\AbSch_k)) \longrightarrow \AbSt_k,$$
It is well known that this
functor is not fully faithful in general, as this can be seen on higher ext-groups. The
ext-groups defined in $\D^b(Pro(\AbSch_k))$ are the \emph{naive extension groups} of commutative
pro-algebraic groups, and can be described in terms of Yoneda extensions. On the other side, 
the extension groups in $\AbSt_k$ are the ext-groups of sheaves. It is know that, when 
$k$ is algebraically closed, 
$Ext^i$ vanish when $i>2$ and when computed in $\D^b(Pro(\AbSch_k))$ (see \cite{MR0266929}), but 
this is not the case in $\AbSt_k$ (see for instance \cite{MR0255550}).

Using the projection $q$ and the $\s$-functor $h$, we get for any commutative group scheme $G$
of finite type over $k$, 
and any integer $i$ a natural morphism of groups
$$\phi^i_G : [J_*(X),G[i]]_{\D^b(Pro(\AbSch_k))} \longrightarrow 
[h(J_*(X)),h(G)[i]] \longrightarrow [\HH^{alg}(X),h(G)[i]],$$
where we have denoted by $[-,-]$ the morphism in the homotopy category 
of stacks
$$[F_1,F_2] = \pi_0(Map_{\St_k}(F_1,F_2)).$$

By universal property the right
hand side can be identified with $H^i_{fppf}(X,G)$. 
We thus get under this restriction
a natural morphism
$$\phi^i_G : [J_*(X),G[i]]_{\D^b(Pro(\AbSch_k))} \to H^i_{fppf}(X,G).$$
For simplicity we denote by $Ext^i_{naive}(J_*(X),G)$ the left hand side of the above
morphism.

\begin{prop}\label{p6}
For all $i\geq 0$ and all commutative
group scheme $G$ which is affine if $i>0$ and unipotent if $i>1$, the morphism
$$\phi_G^i : Ext^i_{naive}(J_*(X),G) \to H^i_{fppf}(X,G)$$
is an isomorphism.
\end{prop}

\textit{Proof.} As in \cite{MR0266929} we can use Galois descent in order to 
reduce to the case where $k$ is algebraically closed.
We start by using that 
the projection $q : \HH^{alg}(X) \to J_*(X)$ is compatible with the filtrations on both sides: 
the filtration by dimension on $\HH^{alg}(X)$ and the naive filtration on $J_*(X)$. To prove the proposition
we will proceed by descending induction on $d$ and will prove that the induced
morphism 
$$F_d(q) : F_d\HH^{alg}(X) \to J_{\geq n-d}(X)$$
induces isomorphisms
$$[J_{\geq n-d}(X),G[i]]_{\D^b(Pro(\AbSch_k))} \simeq [F_d\HH^{alg}(X),G[i]]$$
for all $d$, and all $G$ and $i$ as in the statement of the proposition. Using the fibration sequences
$$F_{d+1}\HH^{alg}(X) \to F_{d}\HH^{alg}(X) \to Gr_{d}\HH^{alg}(X) \qquad J_{\geq n-d-1}(X) \to J_{\geq n-d}(X) \to 
J_{d}(X)[d]$$
and the long exact sequences obtained after applying $[-,G[i]]$, we see that by induction on $d$
we are thus reduced to show that, the induced morphism
$$Ext^{i}_{naive}(J_d(X)[d],G) \to [Gr_{d}\HH^{alg}(X),G[i]]$$
is an isomorphism (for all $d$, and for all $G$ and $i$ as in the statement of the proposition). 
By the proposition \ref{p4}, the right hand side can be identified with 
local cohomology $\bigoplus_{x\in X^{(n-j)}}H^{d+i}_x(X,G)$. We also note that 
by construction $J_d(X)=\prod_{x\in X^{(n-d)}}J_x(X)$, and thus we are reduced to show that 
for a given  point $x \in X$, of codimension $d$, and all $i\geq 0$, the natural morphism
$$\psi^i_G : Ext^i_{naive}(J_x(X),G) \to H^{d+i}_x(X,G),$$
is bijective
under the extra condition stated in the proposition: $G$ is affine if $d+i>0$ and unipotent if $d+i>1$.
Recall here that $J_x(X)$ corepresents the functor $G \mapsto H^{d}_x(X,G)$, 
where $d$ is the codimension of $x$ in $X$.

We first deal with the two special cases $i+d=0$ and $i+d=1$. When $i+d=0$, and thus
$i=d=0$, we have $Ext^0(J_x(X),G)=Hom(J_x(X),G) 
\simeq H^d_x(X,G)$ and it is then obvious that $\psi^0_G$ is bijective. Suppose now that $d+i=1$, 
and thus $G$ is assumed to be affine. If $x$ is of codimension $1$, we have $i=0$, and 
again $\psi_G^0$ is bijective by definition of $J_x(X)$. If $x$ is of codimension $0$
and $i=1$, we have to prove that the natural morphism
$$Ext^1_{naive}(J_0(X),G) \to H^{1}_{fppf}(K(X),G),$$
where $K(X)$ is the fraction field of $X$. To show this we start by the following lemma
treating the cases of $G=\Ga$ and $G=\Gm$.

\begin{lem}\label{lp6}
With the notation above the morphism
$$\psi_G^1 : Ext^1_{naive}(J_0(X),\Ga) \to Ext^1_{naive}(J_0(X),\Gm)\simeq 0$$
is bijective.
\end{lem}

\textit{Proof of the lemma.} Suppose that we are given an extension of commutative 
pro-algebraic groups
$\xymatrix{0 \ar[r] & \Ga \ar[r] & H \ar[r] & J_0(X) \ar[r] & 0.}$ By the universal 
property of $J_0(X)$, constructing a section $J_0(X) \to H$ is equivalent to give 
a point $s \in H(K(X))$, which projects to the universal element in $J_0(X)(K(X))$. 
But we have an exact sequence in $fppf$ cohomology of $K(X)$
$H(K(X)) \to J_0(X)(K(X)) \to H^1_{fppf}(K(X),\Ga)=0$, showing that $H(K(X)) \to J_0(X)(K(X))$
is surjective. Therefore the short exact sequence splits showing that $Ext^1_{naive}(J_0(X),\Ga)=0$.
The proof for $\Gm$ is similar.
\hfill $\Box$ \\

For a general commutative affine group scheme $G$ over $k$, we can decompose 
$G$ as a direct sum $G^u\oplus G^m$, of a unipotent part and a multiplicative part. 
We can moreover use compositions series on each components $G^u$ and $G^m$ to reduce to the case
where $G$ is either $\Ga$ or $\Gm$, or a subgroup of $\Ga$ or $\Gm$ with short
exact sequences
$$0 \to G \to \Gm \to \Gm \to 0 \qquad 0 \to G \to \Ga \to \Ga \to 0.$$
We then consider the action of $\psi_G^1$ on the
associated long exact sequences, and the lemma \ref{lp6} shows that indeed
$\psi_G^1$ is bijective for any $G$. 

We now turn to the case where $i+d>1$ and $G$ is unipotent, which is proven in a very similar fashion
than the case $i=1$ and $d=0$. To begin with, we prove that $\psi_G^i$ is an isomorphism for
$G=\Ga$. Indeed, when $i=0$ this is obvious by the definition of 
$J_x(X)$. For $i>0$ the domain and codomain of $\psi_{\Ga}^i$ vanish. Indeed, 
$H^{d+i}_x(X,\Ga)\simeq H^{d+i}_x(X,\OO)\simeq 0$ for $i>0$ by local duality
for coherent sheaves. On the other hand, $Ext^i_{naive}(J_x(X),\Ga)=0$ by the following lemma.

\begin{lem}\label{lp62}
With the notation above and for all $i>0$
$$\psi_{\Ga}^i : Ext^i_{naive}(J_x(X),\Ga) \simeq 0.$$
\end{lem}

\textit{Proof of the lemma.} Suppose that we are given an extension of commutative 
pro-algebraic groups
$\xymatrix{0 \ar[r] & \Ga \ar[r] & H \ar[r] & J_x(X) \ar[r] & 0.}$ By the universal 
property of $J_0(X)$, construction a section $J_0(X) \to H$ is equivalent to give 
an element $s \in H^{d}_x(X,H)$, which projects to the universal element in $H^d_x(X,J_x(X))$. 
But we have an exact sequence in $fppf$ cohomology of $X$ local at $x$
$H^{d}_x(X,H) \to H^d_x(X,J_x(X)) \to H^{d+1}_x(X,\Ga)=0$, showing that $H^{d}_x(X,H) \to H^d_x(X,J_x(X))$
is surjective. Therefore the short exact sequence splits showing that $Ext^1_{naive}(J_x(X),\Ga)=0$.
Moreover, it is well known, and can be proven using Dieudonné theory, that 
$Ext^i_{naive}(A,\Ga)=0$ for all $i>1$ and all commutative unipotent pro-algebraic group $A$ over $k$
(see \cite{MR0302656} proposition V-1 5.1 and 5.2).
\hfill $\Box$ \\

We thus know that $\psi_{\Ga}^i$ is an isomorphism for all $i$. It remains to treat
the case of a general commutative unipotent group $G$.
Any such group possesses a composition series whose layers
$G_\alpha$ can be realized in a short exact sequence 
$0 \to G_\alpha \to \Ga\to \Ga$. The associated long exact sequences and the lemma
easily imply that $\psi_G^i$ is also an isomorphism for all $i$ 
(note that the domains and codomains of $\psi_G^i$ can be non-zero for $i=1$).
\hfill $\Box$ \\

\begin{rmk}
\begin{enumerate}
    \item \emph{The 
    previous proposition provides a close relation between our algebraic homology $\HH^{alg}(X)$
and the original Grothendieck's construction $J_*(X)$. The latter is a truncated version
of the former, but both can be used in order to express the fppf cohomology of $X$ 
with coefficients in a commutative group scheme. }
\item \emph{The extra conditions on $G$, being affine for $i>0$ and unipotent for $i>1$, do not
appear in Grothendieck's original construction. This is 
because his construction is originaly made using the Zariski topology instead ot he fppf
topology, and higher Zariski cohomology groups with coefficients in group schemes vanish outside
of the unipotent context.}

\end{enumerate}

\end{rmk}

\section{Partial generalization to more general schemes}

In this section we partially extend the results of the previous sections to
the more general case of schemes over arbitrary bases. For this, we restrict ourselves
to the unipotent part of algebraic homology and use the formalism of commutative groups 
in affine stacks (in the sense of \cite{chaff}) as replacement for 
(pro) commutative group stacks. On the way, we develop the general framework of \emph{affine homology}, 
the natural homology theory associated to affine stacks. \\

We fix a base commutative ring $k$. We will work from now with the fpqc topology on affine $k$-schemes, 
so $\St_k$ will denote the $\s$-category of hypercomplete fpqc-stacks on $\Aff_k$.

\subsection{Affine homology}

When $C=\St_k$ is the $\s$-category of stacks
over $k$, the $\s$-category of abelian group objects in $C$ (also called
\emph{group-complete strictly commutative monoids}, see \cite[\S 7]{MR4568994}) 
is naturally equivalent to $\AbSt_k$ (already used previously in this work, but
in the setting of the fppf topology), and is thus
identified with the $\s$-category of connective complexes of sheaves of abelian groups on the 
site of affine $k$-schemes with the fpqc topology.

\begin{df}\label{d7}
The $\s$-category of \emph{commutative affine group stacks over $k$} is defined to be
the full sub-$\s$-category of $\AbSt_k$ consisting of all commutative group stacks $E$ 
such that for all $n\geq 0$ the underlying stack $|E[n]| \in \St_k$ is an affine stack
in the sense of \cite{chaff}. The $\s$-category of all commutative affine group stacks
will be denoted by $\AbChAff_k$.
\end{df}

We caution the reader that the notation $\AbChAff_k$ can be misleading. Indeed, 
the $\s$-category $\ChAff_k$ defines another $\s$-category $Ab(\ChAff_k)$, of 
abelian group objects in $\ChAff_k$. Recall that this is formally defined as the full sub-$\s$-category 
of $Fun^{\s}(T^{op},\ChAff_k)$, consisting of all $\s$-functors sending finite sums in $T$ to finite products in $C$, 
and where $T$ is the algebraic theory of abelian group (i.e. the opposite of the
category of free abelian group of finite rank). The $\s$-category $\AbChAff_k$ is smaller
than $Ab(\ChAff_k)$. Indeed, any affine commutative group scheme $G$ defines an object 
in $Ab(\ChAff_k)$. However, we will see below (see Proposition \ref{p7}) 
that $\Gm$ is not an object in $\AbChAff_k$. 
We thus have canonical fully faithful $\s$-functors
$\AbChAff_k \hookrightarrow Ab(\ChAff_k) \hookrightarrow \AbSt_k$, but these three $\s$-categories are all
distinct. 

Note also that because affine stacks are stable by small limits in $\St_k$, 
$\AbChAff_k \subset \AbSt_k$ is clearly stable by small limits. \\

An important first observation is the following statement, which identifies
the $\s$-category $\AbChAff_k$ when $k$ is a field. It is the commutative group
version of the characterization of affine stacks by means of their
homotopy sheaves which is proven in \cite{chaff}.

\begin{prop}\label{p7}
If $k$ is a field, then the natural inclusion
$$\AbChAff_k \subset \AbSt_k$$
identifies $\AbChAff_k$ with the sub-$\s$-category 
consisting of all $E$ such that for all $i\geq 0$ the sheaf $\pi_i(E)$ is representable by 
a commutative and unipotent affine group scheme over $k$.
\end{prop}

\textit{Proof.} We start by the following lemma. 

\begin{lem}\label{lp7}
Let $E$ be an affine stack over $k$ which is endowed with a group structure. Then, 
for all $i>0$ the sheaves $\pi_i(E)$ are representable by affine unipotent group schemes. 
\end{lem}

\textit{Proof of the lemma.} This is already contained in the proof of \cite[Cor. 3.2.7]{chaff}
which we reproduce here for simplicity. Let $A=\OO(E)$ be the commutative cosimplicial algebra
corresponding to $E$, which is endowed with a $H_\s$-Hopf algebra structure in the
sense of \cite[Def. 3.1.4]{chaff}, or in other words which is 
a cogroup in the $\s$-category of commutative cosimplicial $k$-algebras. As a consequence,
$H^*(A)$ is a graded commutative Hopf $k$-algebra and therefore each $H^0(A)$-module $H^i(A)$
is free (because it corresponds to an equivariant quasi-coherent sheaf
on the group scheme $\Spec\, H^0(A)$), and in particular $A$ is flat as a complex of $H^0(A)$-modules.
Let $G=\Spec\, H^0(A)$ and $p : E \to G$ the natural projection induced
by the canonical morphism $H^0(A) \to A$. The fiber of $p$ is $\Spec\, A'$, where 
$A'=A \otimes_{H^0(A)}k$, which, by the flatness of $A$ over $H^0(A)$, is reduced $H^0(A')\simeq k$.
By \cite[Thm. 2.4.5]{chaff} the corresponding affine stack $E'=\Spec\, A'$ is 
pointed and connected, and thus all the sheaves
$\pi_i(E')$ are representable by unipotent group schemes. The long exact sequence in homotopy therefore implies 
the lemma, as $\pi_i(E)\simeq \pi_i(E')$ for all $i>0$.
\hfill $\Box$ \\

The lemma \ref{lp7} has the following important consequence: for any 
object $E \in \AbChAff_k$, the homotopy groups $\pi_i(E)$ are all representable by unipotent group
schemes, for any $i\geq 0$. Indeed, by the lemma we already know that this is the case for 
$i>0$. We now consider $E[1]$,  and by definition we know that $E[1]$ is again a group
object in affine stacks. In particular the lemma implies that $\pi_1(E[1]) \simeq \pi_0(E)$
is representable by a unipotent group scheme.

Conversely, let $E \in \AbSt_k$ be an object such that the sheaves $\pi_i(E)$
are representable by unipotent group schemes. We have to show that 
each stack $|E[n]|\in \St_k$ is affine. But, the classifying stack
$B(|E[n]|)$
is affine because of \cite[Thm. 2.4.1]{chaff}, and thus so is $\Omega_*B(|E[n]|) \simeq
|E[n]|$ as affine stacks are stable by limits.
\hfill $\Box$ \\

As a direct consequence of the proposition \ref{p7} we have the following corollary, which is the
abelian analogue of the characterizations of affine stacks as being limits of stacks of the form
$K(\Ga,n)$.

\begin{cor}\label{cp7}
We continue to assume that $k$ is a field.
The sub-$\s$-category $\AbChAff_k \subset \AbSt_k$ is the smallest sub-$\s$-category which is
stable by small limits and contains all the objects $\Ga[n]$ for all $n\geq 0$.
\end{cor}

\textit{Proof of the corollary.} Let $C \subset \AbSt_k$ be the smallest 
sub-$\s$-category containing $\Ga[n]$ and stable by small limits.
Clearly $C$ is contained in $\AbChAff_k$, and we know that $\AbChAff_k$ is stable by limits and
that it contains all the objects $\Ga[n]$. The converse is
a direct consequence of Postnikov decomposition and the proposition. By 
\cite[Prop. B1]{MR4191417}, Postnikov towers always converge for the fpqc topology (see also 
\cite{mondal2022postnikov} for a more general statement), and thus any 
$E \in \AbChAff_k$ can be written as $E \simeq \lim_n \tau_{\leq n}E$. 
Moreover, the proposition implies that each $\tau_{\leq n}(E)$ belongs to $\AbChAff_k$. Finally, 
we proceed by induction on $n$ to show that each $\tau_{\leq n}(E)$ belongs to $C$, by means of the cartesian squares
$$\xymatrix{
\tau_{\leq n}(E) \ar[r] \ar[d] & 0 \ar[d] \\
\tau_{\leq n-1}(E) \ar[r] & \pi_{\leq n}(E)[n+1].
}$$
As $\pi_{\leq n}(E)$ is a unipotent group scheme by the proposition \ref{p7}, it is itself contains 
in $C$, and thus so is $\pi_{\leq n}(E)[n+1]$. \hfill $\Box$ \\

We do not know if the corollary \ref{p7} remains valid over an arbitrary base ring $k$. We thus
restrict our notion of commutative affine group stacks to a notion of \emph{affine homology types}, 
by considering all the objects generated by limits of $\Ga[n]$.

\begin{df}\label{d7'}
The \emph{$\s$-category of affine homology types over $k$} is the smallest 
sub-$\s$-category of $\AbSt_k$ which is stable by small limits and 
contains the objects $\Ga[n]$ for all $n\geq 0$. It is denoted by
$$\widehat{\AbChAff_k} \subset \AbSt_k.$$
\end{df}

As a first observation, all the objects $\Ga[n]$ are commutative affine group stacks, and thus
we clearly have $\widehat{\AbChAff_k} \subset \AbChAff_k$. However, outside of the case  of a base
field we do not know if this inclusion is strict. In any case, we will see below that 
$\widehat{\AbChAff_k}$ is better behaved than $\AbChAff_k$, thanks to the following proposition. 
We do not know if the similar statement holds true for $\AbChAff_k$.

\begin{prop}\label{p8}
The
forgetful $\s$-functor of $\s$-categories $\widehat{\AbChAff_k} \hookrightarrow \ChAff_k$,
sending $E$ to $|E|$, possesses a left adjoint. 
\end{prop}

\textit{Proof.} We consider the inclusion of opposite $\s$-categories. 
We know that $(\ChAff_k)^{op}$ is equivalent to the $\s$-category associated with the 
combinatorial model category of commutative cosimplicial $k$-algebras by means of the $\Spec$ functor (see 
\cite{chaff}). Therefore, $(\ChAff_k)^{op}$ is a presentable $\s$-category. Moreover, $(\widehat{\AbChAff_k})^{op}$
is generated by colimits by the small set of objects $\Ga[n]$. Moreover, when considered as objects
in $(\widehat{\AbChAff_k})^{op}$, the objects $\Ga[n]$ are small, or more precisely 
$\Ga[n]$ are cosmall objects in $\AbChAff_k$ as this is implied by the following lemma.

\begin{lem}\label{lp8}
Let $C$ be a complete $\s$-category and $Ab(C)$ be the $\s$-category of 
abelian group objects in $C$. Let $G \in Ab(C)$ be an object satisfying the following two 
conditions.

\begin{enumerate}
    \item The underlying object $G$ is cocompact in $C$: $Map_C(-,G)$ sends filtered limits to 
    filtered colimits.
    \item The object $G$ is truncated in $C$: there is an $n$ such that $Map(E,G)$ is $n$-truncated
    for all $E \in C$.
\end{enumerate}
Then, $G$ is cocompact as an object in $Ab(C)$.
\end{lem}

\textit{Proof of the lemma.} The mapping spaces between two abelian group objects $H$ and $G$ in $C$ can be
computed as a certain homotopy limit over $\Delta$ of the form
$$\Map_{Ab(C)}(H,G) \simeq \lim_{i\in \Delta} \Map_C(H^{n(i)},G^{m(i)}),$$
where $i \mapsto n(i)$ and $i \mapsto m(i)$ are certain functions. The existence of such 
a limit decomposition can be seen for instance using MacLane resolutions of abelian group objects 
as presented in \cite[\S 2]{MR0258842}.
Because $G$ is $n$-truncated this homotopy limit reduces to a finite homotopy limit, which in turn
will commute with filtered colimit. The result follows easily as all objects $G^{m(i)}$ are cocompact in $C$.
\hfill $\Box$ \\

The lemma implies that $(\widehat{\AbChAff_k})^{op}$ is a presentable $\s$-category 
It is therefore a formal consequence of Freyd representability theorem (see \cite[Thm. 4.1.3]{MR4093970}) that 
the forgetful $\s$-functor $(\widehat{\AbChAff_k})^{op} \rightarrow (\ChAff_k)^{op}$, which commutes
with small colimits, admits a right adjoint. \hfill $\Box$ \\

We warn the reader that the forgetful $\s$-functor $\AbSt_k \rightarrow \St_k$
also possesses a left adjoint, given by $F \mapsto \ZZ \otimes F$, but
that the square 
$$\xymatrix{
\ChAff_k \ar[r] \ar[d] & \St_k \ar[d] \\
\widehat{\AbChAff_k} \ar[r] & \AbSt_k}$$
of course does not commute, even when $k$ is a field. 
Indeed, the construction $\ChAff_k \to \widehat{\AbChAff_k}$ involves
considering some free abelian group objects in $\ChAff_k$, and thus involves certain colimits, 
which are computed in $\ChAff_k$ (and the inclusion $\ChAff_k \subset \AbSt_k$ does not
preserve colimits).

We also remind that there are set theory issues in the situation as the site of affine schemes
$Aff_k$ is not a small site. As a result, the inclusion $\ChAff_k \subset \St_k$ does not 
possess a left adjoint as objects in $\St_k$ might not be small. However, 
we can restrict ourselves to small stacks $\St_k^{sm} \subset \St_k$, which are all
objects that can be written as small colimits of affine schemes. By definition 
affine stacks are small, as they are spectrum of small commutative cosimplicial $k$-algebras, so 
we have an inclusion $\ChAff_k \hookrightarrow \St_k^{sm}$. This last inclusion 
does have a left adjoint, called the \emph{affinization} and denoted by $F \mapsto (F\otimes k)^u$. 
It is explicitely given by considering $\OO(F)$, the commutative cosimplicial $k$-algebra of functions
on $F$, which is small as $F$ is assumed to be small, and setting $(F\otimes k)^u \simeq \Spec\, \OO(F)$.

As a result, and including our proposition \ref{p8}, we get an $\s$-functor 
$$\HH^u(-) : \St_k^{sm} \longrightarrow \widehat{\AbChAff_k},$$
which is a left adjoint to the forgetful $\s$-functor $|-| : \widehat{\AbChAff_k} \rightarrow \St_k^{sm}$.

\begin{df}\label{d8}
For a small stack $X \in \St_k$, its \emph{affine homology} is defined 
as $\HH^u(X)$.
\end{df}

Affine homology can be hard to compute. We give here the simplest example
namely the affine homology of a constant stack with finite $1$-connected fibers. 
Let $K$ be a $1$-connected finite simplicial set, which is considered as
a constant stack $\underline{K} \in \St_k$. As such it is a
colimit of the punctual stack $*$ taken over the category of simplices in $K$. Therefore, 
its affine homology is given as
$$\HH^u(\underline{K}) \simeq \underset{\Delta(K)}{colim}\HH^u(*),$$
and we are thus reduced to compute the affine homology of the point. For this, we recall
the existence of the commutative group scheme $\HH$ of \cite{MR4191417}, which can be
defined as the spectrum of the Hopf algebra of integral valued polynomials. Equivalently, 
$\HH$ sits inside the group scheme of big Witt vectors $\WW$ as the sub-functor of
simultaneously fixed points by all Frobenius $F_p$. The group scheme
$\HH$ is the additive part of ring scheme structure, for which the multiplication is
induced by the multiplication of Witt vectors. 

We claim that $\HH^u(*)\simeq \HH$, considered as a commutative group for the additive
structure. First of all, as shown in \cite{MR4191417}
all the $K(\HH,n)$ are abelian group objects in affine stacks, and thus
$\HH$ is indeed an object of $\AbChAff_k$. Using the presentation as fixed points
of Frobenius on $\WW$, it can be checked that $\HH$ lies furthermore in $\widehat{\AbChAff_k}$
and thus is an affine homology type.
In \cite{MR4191417} is also known that when $n>0$
$K(\HH,n)$ is the affinization of the space $K(\ZZ,n)$, and thus
the natural morphism $K(\ZZ,n) \to K(\HH,n)$, which is a morphism of abelian group objects
induced by the canonical map $\ZZ \to \HH$, makes $K(\HH,n)$ into the affinization 
of $K(\ZZ,n)$. Indeed, for any abelian group object $E$ in $\ChAff_k$, 
the mapping space $Map_{\AbSt_k}(K(\ZZ,n),E)$ 
can be computed by a certain homotopy end 
$$Map_{\AbSt_k}(K(\ZZ,n),E) \simeq hoend \left((i,j) \mapsto  Map_{\St_k}(K(\ZZ,n)^i,E^j)\right),$$
where $i$ and $j$ lives in the category of free abelian groups of finite rank. 
Similarly, we have
$$Map_{\AbChAff_k}(K(\HH,n),E) \simeq hoend \left((i,j) \mapsto  Map_{\ChAff_k}(K(\HH,n)^i,E^j)\right).$$
The affinization of $K(\ZZ,n)^i\simeq K(\ZZ^i,n)$ being $K(\HH^i,n)$ (see \cite{MR4191417}), 
we deduce easily that the natural morphism
$$Map_{\AbChAff_k}(K(\HH,n),E) \longrightarrow Map_{\AbSt_k}(K(\ZZ,n),E)\simeq \Omega_*^n(E(k))$$
is an equivalence for any $E \in \AbChAff_k$. In other words, the reduced affine homology 
of $S^n$ is $K(\HH,n)$, for any $n>0$. But we then have 
$Map_{\AbChAff_k}(\HH,E)\simeq Map_{\AbChAff_k}(K(\HH,1),K(E,1))\simeq \Omega_*(K(E,1)(k))\simeq E(k)$, and thus
$\HH$ is indeed the affine homology of $*$ as wanted.

To put things in a more concise manner, we can write, with the notations and conditions above
$$\HH^u(\underline{K}) \simeq K\otimes \HH,$$
where the tensored structure is computed in the $\s$-category $\AbChAff_k$ of commutative affine group stacks over
$k$. \\

For our purpose, two general properties of affine homology are of some importance, 
namely base change and filtration by dimension. \\

\textbf{Base change.} As a start, affine homology is base sensitive, and 
should rather be written as $\HH^u(X/k)$ in order to mention the base ring.
Let $k \to k'$ be a morphism of commutative rings, 
and $X \in \St_k^{sm}$ a small stack over $k$. We denote by $\otimes_k k'$ the base
change from stacks over $k$ to stacks over $k'$, and similarly for affine stacks. 
The base change commutes with limits and thus induces a well defined base change
$\AbChAff_k \to \AbChAff_{k'}$. As this sends $\Ga$ to $\Ga$, base change restricts to 
an $\s$-functor
$$\otimes_k k' : \widehat{\AbChAff_k} \to \widehat{\AbChAff_{k'}}.$$

Therefore, by universal property we have a canonical morphism of 
affine homology types over $k'$
$$\HH^u(X\otimes_k k'/k') \to \HH^u(X/k)\otimes_k k'.$$
This morphism is not an equivalence in general, but we claim it is so 
when $k'$ is finite and of local complete intersection over $k$. Indeed, by the very definition 
of affine homology this would follow formally from the
fact that the direct image $\s$-functor $R_{k'/k} : \AbSt_{k'} \to \AbSt_k$, right adjoint
to the base change, preserves affine homology types. As these are generated by limits from the
$\Ga[n]$, it is enough to show that $R_{k'/k}(\Ga[n])$ is an affine homology type 
over $k$. But $R_{k'/k}(\Ga[n])$ is nothing else than
$|p_*(\Ga)[n]|$, where $p : \Spec\, k' \to \Spec\, k$. Moreover,
by assumption on $k'/k$ the stack $p_*(\Ga)[n]$ is a linear stack over $k$ 
associated to the perfect complex $E:=(k')^\vee[-n]$, which is the dual 
of $k'[n]$. More explicitly, $p_*(\Ga)[n]$ can be written as the spectrum
of the cosimplicial commutative $k$-algebras $Sym_{k}(E)$, with $E$ 
a perfect complex of $k$-modules. These linear stacks, as objects in $\AbSt_k$
can be decomposed by considering a projective resolution of $E$, 
and are easily seen to be expressible by limits of stacks of the form $V[n]$
where $V$ is a vector bundle on $\Spec\, k$, which are retracts of $\Ga^r[n]$
and thus are affine homology types over $k$.

An important case of base change is when $k$ is regular, and $k \to k(p)$ is the 
quotient corresponding to a closed point $p \in \Spec\, k$. In this case base change
tells us that $\HH^u(X/k)$ specializes to $\HH^u(X_p/k(p))$ when pull-backed to the
point $p$. \\

\textbf{Filtration by dimension.} We let $X$ be a $k$-scheme separated and of finite type
over $k$. For simplicity we assume that $k$ is noetherian and $X$ is of finite dimension $d$.
We can proceed as in the definition \ref{d5} and construct a finite filtration on $\HH^u(X/k)$
$$\xymatrix{
0=F_{d+1}\HH^u(X/k) \ar[r] & F_d\HH^u(X/k) \ar[r] & \dots \ar[r] & F_1\HH^u(X/k) \ar[r]
& F_0\HH^u(X/k) = \HH^u(X/k).}$$
The associated graded for this filtration, $Gr_d\HH^u(X/k)$, defined considering cofibers
in $\widehat{\AbChAff_k}$, are affine homology types over $k$ 
which are related to local cohomology in dimension $d$ as follows. We first 
introduce affine homology types of finite presentation, as being 
the objects $E \in \widehat{\AbChAff_k}$ being obtained 
by finite limits and retracts of objects of the form $\Ga[n]$. As we have already seen in 
the lemma \ref{lp8} the objects $\Ga[n]$ are cocompact in $\AbChAff_k$, and thus
any affine homology type of finite presentation are cocompact objects among
affine homology types.

Let then $E$ be an affine homology type of finite presentation. We then have
$$\Map_{\AbChAff_k}(Gr_d\HH^u(X),E)) \simeq \oplus_{x \in X^{(d)}}\HH_x(X,E).$$
Because affine homology type of finite presentation
generate the whole $\s$-category $\widehat{\AbChAff_k}$ by limits, we deduce that
we have a general formula in terms of local affine homology
$$Gr_d\HH^u(X) \simeq \prod_{x\in X^{(d)}}\HH^{u}_x(X),$$
where each factor $\HH^{u}_x(X)$ is defined as the cofiber, taken in 
$\widehat{\AbChAff_k}$ of the natural morphism
$\HH^u(X_x-\{x\}) \to \HH^u(X_x)$.

\subsection{Affine homology as a dg-module}

We consider the object $\Ga \in \AbSt_k$ and the corresponding $\ZZ$-linear dg-algebra of endomorphisms
$$C_{\Ga}:=\mathbb{R}End(\Ga).$$
We get this way the usual adjunction of $\s$-categories
$$M : (\AbSt_k)^{op} \leftrightarrows C_{\Ga}-\dg : D,$$
where $M(E):=\mathbb{R}Hom(E,\Ga)$. As $\Ga[n]$ generates $\cAbChAff_k$ by limits, and as
$\Ga[n]$ are cocompact objects, it is formal to see that for any $E \in \cAbChAff_k$ the adjunction morphism
$$E \longrightarrow D(M(E))$$
is an equivalence. As a result we get the following corollary.

\begin{cor}
The $\s$-functor $M$ restricts to a full embedding
$$M : (\cAbChAff_k)^{op} \longrightarrow (C_{\Ga}-\dg)^{op}.$$
Its essential image consists of all $C_{\Ga}$-dg-modules which are colimits of 
$C_{\Ga}$.
\end{cor}

The above corollary implies that for a small stack $X$ its affine homology $\HH^u(X/k)$ is characterized
by $\rg(X,\Ga)$, as a dg-module over $C_{\Ga}$. The dg-algebra $C_{\Ga}$ is however quite complicated in general. 
It simplifies when $\Ga$ is restricted to sites of perfect schemes (see \cite{MR0638602}). In particular, 
things tend to simplify when the base ring $k$ is perfect (or at least perfectoid). Another simplification
consists of using the group of big Witt vectors $\WW$ instead of the additive group $\Ga$, as for instance
$Ext^1(\WW,\WW)=0$, and again higher ext vanish when computed over the site of perfect schemes. Moreover, 
$End(\WW)$ is the usual Dieudonné ring so dg-modules over $\mathbb{R}End(\WW)$ are closed related to 
usual Dieudonné modules. However, 
$\WW$ is not a cocompact object in $\cAbChAff_k$, and therefore the corresponding adjunction 
with dg-modules is not as well behaved and relates dg-modules over $\mathbb{R}End(\WW)$ and certain 
pro-objects inside $\AbChAff_k$.

\bibliographystyle{plain}
\bibliography{Biblio.bib}

\noindent Bertrand To\"{e}n, {\sc IMT, CNRS, Universit\'e de Toulouse, Toulouse (France)}\\
Bertrand.Toen@math.univ-toulouse.fr \\

\end{document}